\documentclass[12pt]{article}
\usepackage{amsmath,amssymb,amsthm}

\usepackage{fullpage}

\title{On Those Analytic Functions of One Argument Which Possess an Algebraic Addition Theorem}

\author{PAUL KOEBE \\[12pt]
        Leipzig\footnote{(translator)The following is a literal word-for-word translation from the German of the  the first five articles of  the work ``Ueber diejenigen analytischen Funktionen eins Arguments, welche ein algebraisches Addititionstheorem besitzen, ...," which appears in the book \emph{Mathematische Abhandlungen, \textsc{Hermann Amandus Schwarz} su seinem fuenfzigjaehrigen Doktorjubilaeum, am 6 August 1914, gewidmet von Freunden und schuelern}, 1974, Chelsea , pp 192-208.}}

\date{\today}

\theoremstyle{definition}

\begin{document}
\maketitle

In his lectures held at the university of Berlin on the theory of elliptic functions, \textsc{Weierstrass} started with the problem: to determine all analytic functions of one argument which possess an algebraic addition theorem.  \textsc{Weierstrass} carried out the complete solution of this problem in his lectures.  To determine the multiple-valued functions with an algebraic addition theorem he used different methods at different times.  Essentially it comes down to showing that the number of branches of such a function is finite.

This paper contains, what I believe to be, the simplest and most natural way to prove that a multiple-valued function with an algebraic addition theorem is an algebraic function of a single-valued function.  It is likely that \textsc{Weierstrass} must have come very close to this method of proof, for the tools which are used in the method of proof are very close to the ones he used for this purpose, something I learned from the lectures of \textsc{H. A. Schwarz}.

In the seventh volume of \emph{Acta mathematica} \textsc{Phragmen} published a method which he found for the solution of the problem of \textsc{Weierstrass}.\footnote{\textsc{Phragmen}: ``Sur un theoreme concernant les fonctions elliptiques".  Acta Matematica, Vol 7.  The extension of the probem of \textsc{Weierstrass} to functions of several variables forms the content of a paper by \textsc{P. Painleve}: ``Sur les fonctions qui admettent un theoreme d'addition".  Acta Mathematica, Vol 27}  Although the method of proof which I found has some unmistakable relations with that of \textsc{Phragmen} there still are some notable differences.  Since the appearance of my thesis, \textsc{Falk} has investigated the functions of one argument with an algebraic addition theorem in the paper ``Ueber die Haupteigenschaften derjenigen Funktionen eines Arguments, welche Additionstheoreme besitzen" (Nova Acta soc. Upsal., Ser. IV, Vol I, Nr. 8, 1907.)

In \S 1 of this paper we will prove the following theorem found by \textsc{Weierstrass}: Every analytic function $\phi(u)$ of the argument $u$, which possesses an algebraic addition theorem, i.e., which has the property that between any three of the of the values of the function corresponding to the three values of the argument $$u, \\ v, \\ u+v,$$ that is between $$\phi(u), \ \ \phi(v),\ \ \phi(u+v)$$there exists an algebraic equation $$G[\phi(u),  \phi(v),\phi(u+v)]=0$$with coefficients independent of $u$ and $v$ is:\begin{itemize}
  \item either an algebraic function of $u$;
  \item or, if $\mu$ denotes a suitable constant, an algebraic function of the exponential function $e^{\mu u}$;
  \item or, if $\omega,\omega'$ denote suitable constants, an algebraic function of the function $\wp(u|\omega,\omega')$.
\end{itemize}(See \textsc{H.A. Schwarz}: Formln und Lehrsaetze zum  Gebrauche der elliptischen Funktionen (p. 1-2), second edition, Berlin, 1993).

The special assumption we make for this proof, that the given algebraic equation holds for an element of the function $\phi(u)$ which exists and is regular in the neighborhood of the point $u=0$, and made for methodological advantages, will be shown in \S 2 to represent no essential limitation to the generality of the solution.  However it is not appropriate when the formally more general assumption holds, to assume for an analytic function $\phi(u)$ that the algebraic equation
$$G[\phi(u),  \phi(v),\phi(u+v)]=0$$holds for suitable choice of the branches $$\phi(u), \ \ \phi(v),\ \ \phi(u+v).$$  This assumptions carries with it an inconvenience in that one makes the implicit assumption that the elements, whose function values we have denoted with $$\phi(u), \ \ \phi(v),\ \ \phi(u+v)$$can be transformed into each other by analytic continuation along a path.  This inconvenience can be eliminated by carrying over to a more general formulation of the theorem of \textsc{Weierstrass}:

\emph{If between any three analytic functions $$\phi(u), \psi(v), \chi(u+v)$$ there subsists an algebraic equation $$F[\phi(u), \psi(v), \chi(u+v)]=0$$with coefficients independent of $u$ and $v$, then the domain of the argument of each of these functions be extended to al finite values, and, to be sure, all are algebraic functions of $u$, or all algebraic functions of one and the same function $e^{\mu u}$, or all algebraic functions of one and the same function $\wp(u|\omega,\omega')$}.

This paper is broken up into the following sections.

\S 1  The determination of all analytic functions of one argument which possess an algebraic addition theorem.

\S 2  Transition to the algebraic functional equation $$F[\phi(u), \psi(v), \chi(u+v)]=0.$$

\S 3  Proof of the theorem that to each single-valued or multiple-valued function $\phi(u)$ with an algebraic addition theorem there corresponds a unique irreducible equation$$G[\phi(u),  \phi(v),\phi(u+v)]=0.$$.

\S 4  On the relation between two analytic functions of the argument $u$ which possess the same algebraic addition theorem.

\S 5  The degree of the irreducible equation $G[\phi(u),  \phi(v),\phi(u+v)]=0$ withy respect to $\phi(u+v), \ \ \phi(u),\ \ \phi(v)$.

\newpage
\section{The determination of all analytic functions of one argument which possess an algebraic addition theorem.}

Let $\phi(u)$ be an analytic function of the argument $u$ which has an element, $\phi_0(u)$, which is regular in the neighborhood of the point $u=0$, and whose function values $\phi_0(u)$, $\phi_0(v)$, and $\phi_0(u+v)$ satisfy an algebraic equation with coefficients independent of $u$ and $v$:\begin{equation}
\label{AT}
G[\phi_0(u),\phi_0(v),\phi_0(u+v)]=0
\end{equation}

We will call radius of the circle of convergence of the power series in $u$, which represents the regular element $\phi_0(u)$, $r$.

If we put $u=v$ in the equation \eqref{AT}, then for $|u|<\frac{r}{2}$ we obtain an algebraic equation with constant coefficients between $\phi_0(u)$ and $\phi_0(2u)$:\begin{equation*}
\label{ }
G[\phi_0(u),\phi_0(u),\phi_0(2u)]=H_1[\phi_0(u),\phi_0(2u)]=0
\end{equation*}If $G[\phi_0(u),\phi_0(v),\phi_0(u+v)]$ is an irreducible polynomial\footnote{(translator) I translate ``ganze rationale Funktion" as ``polynomial".} in the quantities $\phi_0(u)$, $\phi_0(v)$, and $\phi_0(u+v)$, which may always be assumed, then it cannot happen that all of the coefficients in the polynomial $$H_1[\phi_0(u),\phi_0(2u)]$$ vanish.  Otherwise the polynomial $G$ would be divisible by $\phi_0(u)-\phi_0(v)$.  More generally, it follows from the equation\begin{equation*}
\label{ }
H_1[\phi_0(u), \phi_0(2u)]=0
\end{equation*}that: if $n$ is an arbitrarily large positive whole number, and if $|u|<r$, then between $\phi_0(u)$ and $\phi_0(\frac{u}{2^n})$ there exists an algebraic equation with constant coefficients \begin{equation}
\label{n}
H_n\left[\phi_0\left(\frac{u}{2^n}\right), \phi_0(u)\right]=0
\end{equation}

In this equation, $\phi_0(\frac{u}{2^n})$ is a function of $u$ which  is one valued and has the character of an entire function in the domain $$|u|<2^nr.$$

Therefore, the equation \eqref{n} shows that starting with the regular element $\phi_0(u)$ of $\phi(u)$, it can be analytically continued throughout the interior of a circle, $K_n$, centered at the origin with radius $2^nr$.  By means of this continuation the function $\phi(u)$ will have only a finite number of branches and the function will possess at each point $u_0$ in the interior of the circle $K_n$ the character of an algebraic function, that is, its values can be represented in the neighborhood of the value $u_0$ by a finite number of power series, which, in general, proceed in powers of a quantity $\sqrt[m]{u-u_0}$ ($m$ a positive integer), in which the number of terms with negative powers is finite.  If $\phi_1(u)$ and $\phi_2(u)$ denote any two of these branches, then they satisfy the follwing equations \begin{equation*}
\label{}
H_n\left[\phi_0\left(\frac{u}{2^n}\right), \phi_1(u)\right]=0, \ \ H_n\left[\phi_0\left(\frac{u}{2^n}\right), \phi_2(u)\right]=0
\end{equation*}from which the uniquely determined quantity   $\phi_0(\frac{u}{2^n})$ may be eliminated.

The number of singular points which results from the analytic continuation of the function $\phi(u)$ within the circle $K_n$ is finite.  For, if it were infinite, there would be a limit point of singular points in which the function $\phi(u)$ would no longer have the character of an algebraic function.  But such a point cannot be eithr in the internior nor on the boundary of the circle $K_n$ since as a consequence of the equation$$H_{n+1}\left[\phi_0\left(\frac{u}{2^{n+1}}\right), \phi(u)\right]=0$$the function $\phi(u)$ can be continued beyond the boundary of the circle $K_n$ without ceasing to have the character of an algebraic function.

Since the number of singular points of the function $\phi(u)$ which are in the circle $K_n$ is finite, as we just showed, there exist infinitely many pairs of mutually perpendicular diameters of this circle which do not go through any singular points (different from $u=0$).  We choose a definite one of these pairs.  Let the points of one of these diameters geometrically represent the values of the quantity $u_1$ and the points of the other diameter the values of the quantity $u_2$.  With this definition of $u_1$ and $u_2$, every quantity $u$ whose absolute value is smaller than $2^nr$ can be uniquely represented in the form $$u=u_1+u_2$$.

For all values of $u$, whose absolute value is smaller than $r$, there is a uniquely determined branch of the function $\phi(u)$, namely the branch $\phi_0(u)$.  The analytic continuations of this branch along the diameters within the circle $K_n$ are unique because, by assumption, neither diameter goes through a singular point.  We will denote these uniquely determined continuations of the branch $\phi_0(u)$ by $\phi_0(u_1)$ and $\phi_0(u_2)$.

Now let $u$ describe an arbitrary curve $L$ from the origin to an interior point $u'=u_1'+u_2'$ of the circle $K_n$, which lies entirely within the given circle and which goes through no singular point.  As long as the absolute value of $u$ does not exceed $r$, according to the fundamental assumption, the equation \begin{equation*}
\label{ }
G[\phi_0(u_1),\phi_0(u_2),\phi_0(u)]=0
\end{equation*}holds true.  This equation remains true for all analytic continuations, and therefore, if $\phi_0(u)$ transforms into $\phi_0(u')$ while being continued along the path $L$, whiche into $\phi_0(u_1')$ and $\phi_0(u_2')$, the following equation holds \begin{equation*}
\label{ }
G[\phi_0(u_1'),\phi_0(u_2'),\phi_0(u')]=0
\end{equation*}

In this equation  $\phi_0(u_1')$ and $\phi_0(u_2')$ are independent of the particular path joining the origin to $u'$.  Therefore, the number $h_n$ of branches which can result from analytic continuation of the function $\phi(u)$ within the circle $K_n$ cannot be greater than the degree, $g$, which is independent of $n$, of the fundamental equation \begin{equation*}
\label{ }
G[\phi_0(u),\phi_0(v),\phi_0(u+v)]=0
\end{equation*}with respect to 
$\phi_0(u+v)$.  This observation holds for every $n$ however large.  Therefore, the  positive numbers $h_1, h_2, h_3, \cdots$ satisfy the inequalities\begin{equation*}
\label{ }
h_1\leq h_2\leq \cdots \leq h_n\leq h_{n+1}\leq \cdots \leq g
\end{equation*}from which we may conclude that there exists a positive whole number $N$ with the property that  \begin{equation*}
\label{ }
h_N=h_{N+1}=h_{N+2}=\cdots \text{to infinity}
\end{equation*}that is, the set of all the $h_N$ branches of the function $\phi(u)$ obtained from the element $\phi_0(u)$ by analytic continuation in the interior of the circle $K_N$ represents the totality of all branches of the function $\phi(u)$.  For the transition to further circles of greater radius the analytic continuation of the function $\phi(u)$ results in a simultaneous continuation of the $h_N$ branches which transform into \emph{each other} when they circulate around any newly occurring branch point.

It was shown earlier that any two branches of the analytically continued function $\phi(u)$ within the circle $K_n$ are connected by an algebraic equation with constant coefficients.  We can now assert that this theorem holds everywhere for the function $\phi(u)$.  The proof results immediately if we replace $n$ by $N$.

The previous considerations have proven the following: 

\emph{The domain of the argument of the function $\phi(u)$ can be extended to all finite values without the function ceasing to be of algebraic character.  The number of branches of the function $\phi(u)$ is finite, and every branch is an algebraic function of every other branch. }

Among the elementary symmetric functions of all of the branches of the function $\phi(u)$ there exist at least \emph{one} function $\psi(u)$ which is not a constant.  The function $\psi(u)$ is a single-valued function of the argument $u$ which has the character of a \emph{rational} function for all finite values of $u$.  For $|u|<r$, $\psi(u)$ is an algebraic function of $\phi_0(u)$ since all branches of the function $\phi(u)$ are algbraic functions of the branch $\phi_0(u)$ for $|u|<r$.  Therefore, the single-valued function $\psi(u)$ also has an \emph{algebraic addition theorem}, and the \emph{function $\phi(u)$} can be considered \emph{as an algbraic function of the single-valued function $\psi(u)$.}

Since the function $\psi(u)$ has the character of a rational function for all finite values, $\psi(u)$ is either a \emph{rational function of $u$}, or a function of $u$ for which the point $u=\infty$ is an essential singular point.  In the latter case a procedure indicated by \textsc{Weierstrass} be further completed.

From the property that the function $\psi(u)$has an algebraic addition theorem it follows through differentiation of th equation of the addition theorem for $\psi(u)$ with respect to $u$ and $v$ and the subsequent elimination of the quantities $\psi(u+v)$ and $\psi'(u+v)$, that the derivative $\psi'(u)$ of the function $\psi(u)$ is an algebraic function of $\psi(u)$, and that all higher derivatives are rationally expressible through $\psi(u)$ and $\psi'(u)$.  Let $k$ be the degree of the algebraic equation between $\psi(u)$ and $\psi'(u)$ with respect to $\psi'(u)$.  Since $u=\infty$ is an essential singularity a function value $b$ can be determined inside an arbitrarily chosen two-dimensional domain such that at least $k+1$ distinct values of the argument exist for which\begin{equation*}
\label{ }
\psi(a_1)=\psi(a_2)=\cdots=\psi(a_{k+1})=b
\end{equation*}Among these $k+1$ values of the argument there must exist at least two, $a_\alpha$ and $a_\beta$ for which not only does the equation $\psi(a_\alpha)=\psi(a_\beta)$ but also the equation $\psi'(a_\alpha)=\psi'(a_\beta)$.  Therefore, at the points $a_\alpha$ and $a_\beta$ all of the values of the higher derivatives also coincide, and consequently, if $h$ be a variable quantity whose values are restricted to those whose absolute value are sufficiently small, the expansions of 
$\psi(a_\alpha+h)$ and $\psi(a_\beta+h)$ coincide term by term.  The quantity $a_\alpha-a_\beta$ is accordingly a period of the argument of the function $\psi(u)$,

If $2\omega$ denotes a primitive period of the argument of the function $\psi(u)$ and $\mu$ the quantity $\dfrac{\pi i}{\omega}$, we can put\begin{equation*}
\label{ }
\psi(u)=\overline{\psi}(e^{\mu u})=\overline{\psi}(z)
\end{equation*}The function $\overline{\psi}(z)$ is a single-valued function of the variable $z=e^{\mu u}$ for which only the points $z=0$ and $z=\infty$ could be essential singularities.  If neither $z=0$ nor $z=\infty$ is an essential singualrity, then $\overline{\psi}(z)$ is a rational function of $z$, i.e., $\psi(u)$ is a \emph{rational function of $e^{\mu u}$}.  Otherwise, by repeating the earlier argument for the function $\psi(u)$, we can conclude the existence of a new period $2\omega'$ whose ratio with the period $2\omega$ is not real.  The function $\psi(u)$ is therefore a \emph{doubly periodic function} and can be expressed \emph{rationally by $\wp(u|\omega,\omega')$ , $\wp'(u|\omega,\omega')$} and algebraically through $\wp(u|\omega,\omega')$ alone.

%%%%%%%%%%%%%%%%%%%%%%%%%%%%%%%%%%%%%%%%%%%%%%%%
\newpage
\section{Transition to the algebraic functional equation \\ $F[\phi(u), \psi(v), \chi(u+v)]=0.$}

The special assumption made at the beginning of the proof of \textsc{Weierstrass'} theorem that the function $\phi(u)$ has an element which exists and is regular in the neighborhood of the point $u=0$ and for that element there exists an algebraic addition theorem, is a limitation which can easily be removed.

Suppose that between  three analytic functions $\phi(u), \psi(v), \chi(u+v)$, each of one variable, there exists an algebraic relation $$F[\phi(u), \psi(v), \chi(u+v)]=0,$$ with coefficients independent of $u$ and $v$.  We assume the functions $\phi(u), \psi(v), \chi(u+v)$ exist in the neighborhood of the (respective) points $u=a$, $v=b$, $u+v=a+b$ and are there represented by the ordinary power series \begin{equation*}
\label{ }
\mathbf{P_1(u-a)}, \ \ \mathbf{P_2(v-b)}, \ \ \mathbf{P_3(u+v-a-b)}.
\end{equation*}If we put \begin{equation*}
\label{ }
u-a=x, \ \ v-b=y, \ \  \text{therefore} \ \ u+v-a-b=x+y
\end{equation*}it follows from $$F[\mathbf{P_1(x)},\mathbf{P_2(y)},\mathbf{P_3(x+y)}]=0,$$if we put $x$, and then $y$ equal to zero, that $\mathbf{P_2(x)}$ and $\mathbf{P_3(x)}$ are algebraic functions of $\mathbf{P_1(x)}$ and that therefore there exists an algebraic equation between $\mathbf{P_1(x)}$, $\mathbf{P_1(y)}$, $\mathbf{P_1(x+y)}$:\begin{equation*}
\label{ }
G[\mathbf{P_1(x)},\mathbf{P_1(y)},\mathbf{P_1(x+y)}]=0
\end{equation*}It is possible that the substitutions of the special values $x=0$ or $y=0$ could lead in special cases to an identity.  This case can always be avoided by a suitable variation of the points $a$, $b$, $a+b$ in the neighborhood where the series expansions hold.

The above remarks suffice to prove the theorem, whose statement was given in the introduction, for the functions $\phi(u)$, $\psi(v)$, $\chi(u+v)$.

Assuming that the function elements represented by the series $\mathbf{P_1(u-a)}$ , $\mathbf{P_2(v-b)}$, $\mathbf{P_3(u+v-a-b)}$ are elements of the \emph{same} analytic function $\phi(u)$, in place of the equation $$F[\phi(u), \psi(v), \chi(u+v)]=0,$$ we now have the equation $$F[\phi(u), \phi(v), \phi(u+v)]=0,$$which says that the function $\phi(u)$ has an algebraic addition theorem.  Thus we have dismissed the special assumption made in \S 1.

It is useful to state the following specialization of the theorem about the functions  $\phi(u), \psi(v), \chi(u+v)$:

\emph{If an analytic function $\phi(u)$ has the property that there exists an algebraic function of $\phi(u)$ and $\phi(v)$ which only depends on the sum $u+v$, then $\phi(u)$ is either an algebraic function of $u$, or of a function $e^{\mu u}$, or of a function $\wp(u|\omega,\omega').$}

This theorem corresponds to the assumption of an algebraic equation $$F[\phi(u), \phi(v), \chi(u+v)]=0,$$ between the functions $\phi(u)$, $\phi(v)$, $\chi(u+v)$.

%%%%%%%%%%%%%%%%%%%%%%%%%%%%%%%%%%%%%%%%%%%%%%%%

\newpage

\section{Proof of the theorem that to every single-valued or multiple-valued analytic function $\phi(u)$ with an algebraic addition theorem there belongs only one unique irreducible equation $G[\phi(u),  \phi(v),\phi(u+v)]=0.$}

Let there be given a function with an algebraic addition theorem.  By the theorem of \textsc{Weierstrass} it is an algebraic function either of $u$,  or of a function $e^{\mu u}$, or of a function $\wp(u|\omega,\omega').$  If $\phi(u)$, $\phi(v)$, $\phi(u+v)$ are any branches whatsoever of this function, then between them there exists a definite irreducible algebraic equation 
$G[\phi(u),  \phi(v),\phi(u+v)]=0.$  It can be shown that this equation is \emph{independent} of the choice of branches in the following way.

Let $u_0$, $v_0$, $u_0+v_0$ be a special system of values of the argument which has the property that all the branches of the function under consideration are regular in the neighborhood of the points of the plane corresponding to the values  $u_0$, $v_0$, $u_0+v_0$.  Then it is possible to simultanously analytically continue the system of function values $\phi(u_0)$, $\phi(v_0)$, $\phi(u_0+v_0)$ along a path into the system $\overline{\phi}(u_0)$, $\overline{\phi}(v_0)$, $\overline{\phi}(u_0+v_0)$, where $\overline{\phi}(u)$, $\overline{\phi}(v)$, $\overline{\phi}(u+v)$ is any arbitrarily chosen system of branches.

Let the variable $u$ describe a closed path from $u_0$, the variable $v$ from $v_0$, along which 
$\phi(u_0)$ transforms into $\overline{\phi}(u_0)$, $\phi(v_0)$ into $\overline{\phi}(v_0)$. If by this 
$\phi(u_0+v_0)$ transforms into $\overline{\overline{\phi}}(u_0+v_0)$, then it remains to transform $\overline{\overline{\phi}}(u_0+v_0)$ back into $\overline{\phi}(u_0+v_0)$. To this end let the variation of the variable $u$ be limited to values of the form $u+kt$, the variation of the variable $v$ to values of the form $v+k't'$, where $k$ and $k'$ are two constants with a non-real ratio, and $t$ and $t'$  independent real variables.  Then the corresponding points in the complex plane are two straight lines going respectively through $u_0$ and $v_0$.  To any arbitrarily chosen value $u+v$ corresponds a completely determinate pair of values $t,t'$, therefore a uniquely determined pair of values $u,v$.  Since only a finite number of singular points are in any finite domain in the complex plane, however large, the constants $k$ and $k'$  can be chosen in infinitely many ways so that no branch point lies on either of the lines for the mobile $u$ or $v$, while $u+v$ moves from $u_0+v_0$ on a closed curve in the finite plane, along which $\overline{\overline{\phi}}(u_0+v_0)$ transforms into $\overline{\phi}(u_0+v_0)$.  If the constants $k$ and $k'$ are chosen in accordance with this condition then the motion of $u$ and $v$ corresponding to the above motion of $u+v$ returns the values  $\overline{\phi}(u_0)$, $\overline{\phi}(v_0)$ to the values $\overline{\phi}(u_0)$, $\overline{\phi}(v_0)$, respectively.

Therefore we have shown that the system of branches $\phi(u)$, $\phi(v)$, $\phi(u+v)$ can be transformed into any other arbitrarily chosen system of branches $\overline{\phi}(u)$, $\overline{\phi}(v)$, $\overline{\phi}(u+v)$ by simultaneous analytic continuation, whence the theorem to be proved follows.\footnote{This theorem has been proven in a different way by \textsc{Phragmen} in the work cited above.}  This theorem can also be expressed as follows:

\emph{To any function $\phi(u)$ with an algebraic addition theorem, also if it be multiple-valued, corresponds a single uniqu irreducible algebraic equation $G[\phi(u),  \phi(v),\phi(u+v)]=0.$}
 
%%%%%%%%%%%%%%%%%%%%%%%%%%%%%%%%%%%%%%%%%%%%%%

\newpage

\section{On the relation between two analytic functions of the argument $u$ which have the same algebraic addition theorem.}

According to the results of \S 3, to every function $\phi(u)$ with an algebraic addition theorem belongs only one unique irreducible algebraic equation  $G[\phi(u),  \phi(v),\phi(u+v)]=0.$  Now one can consider the question as to, conversely, how uniquely is the function $\phi(u)$ determined by this equation.

The assumption of irreducibility of the equation $$G[\phi(u),  \phi(v),\phi(u+v)]=0,$$ has the consequence that all systems of values $x$, $y$, $z$ satisfying the equation $G(x,y,z)=0$ are individually transformable by analytic continuation into each other\footnote{see the appendix to my thesis (Berlin 1905)}.  Therefore, if $\psi(u)$ is any other analytic function of $u$ which also satisfies the equation $$G[\psi(u),  \psi(v),\psi(u+v)]=0,$$  Then three related variables, $u_1$, $v_1$, $u_1+v_1$ dependent on $u$, $v$, $u+v$ can be determined in accordance with the equations \begin{equation}
\label{*}
\phi(u)=\psi(u_1), \ \ \phi(v)=\psi(v_1), \ \ \phi(u+v)=\psi(u_1+v_1)
\end{equation}Solving these equations for  $u$, $v$, $u+v$ these equations take the form\begin{equation*}
\label{ }
u=\chi_1(u_1). \ \  v=\chi_2(u_2), \ \ u+v=\chi_3(u_1+v_1)
\end{equation*}so that \begin{equation*}
\label{ }
\chi_1(u_1)+\chi_2(u_2)=\chi_3(u_1+v_1)
\end{equation*}It follows by differentiating with respect to $u_1$ and $v_1$ that\begin{equation*}
\label{ }
\chi_1'(u_1)=\chi_2'(v_1)=\chi_3'(u_1+v_1)
\end{equation*}Since $u_1$ and $v_1$ are independent of one another, the common value of these quantities is a constant $\alpha$.  Therefore, if $a_1$ and $a_2$ denote two new constants, and if we put $u+v=w$ and $u_1+v_1=w_1$, we obtain for $u$, $v$, $w$ expressions of the form 
\begin{equation*}
\label{ }
u=\alpha u_1+\alpha_1, \ \ v=\alpha v_1+a_2, \ \ w=\alpha w_1+a_1+a_2
\end{equation*}If we substitute these expressions into the equation \eqref{*}, these latter take the form \begin{equation*}
\label{ }
\phi(\alpha u_1+\alpha_1)=\psi(u_1), \ \ \phi(\alpha v_1+a_2)=\psi(v_1), \ \ \phi(\alpha w_1+a_1+a_2)=\psi(u_1+v_1)
\end{equation*}In these equations we can replace the quantities $u_1$, $v_1$, $w_1$ by $u$.  From the last two it results that $a_1$ is a period of the function $\phi(u)$, and consequently the first can be put into the form\begin{equation*}
\label{ }
\phi(\alpha u)=\psi(u).
\end{equation*}Therefore we have proven the theorem:

\emph{If $\phi(u)$is a function with an algebraic addition theorem, then all analytic functions with the same addition theorem can be represented in the form $\phi(\alpha u)$ where $\alpha$ is a constant,.}

If the origin of the complex plane be changed, this corresponds to the transformation of the function $\phi(u)$ into the function $\phi(\alpha+u)$, and so we can answer the question of when two functions $\phi(\alpha+u)$ and $\phi(\beta+u)$ possess the same addition theorem.  The previous theorem offers us the answer to this question:

\emph{Two functions $\phi(a+u)$ and $\phi(b+u)$ have the same algebraic addition theorem if and only if by a linear transformation $u'=\alpha u+\beta$ of the variable $u$, by which the $a$ is transformed into $b$, the function $\phi(u)$ is transformed into itself.} 

That is, if it is possible to determine two constants $\alpha$ and $\beta$ that the following two equations hold simultaneously:\begin{equation*}
\label{ }
\phi(\alpha u+\beta)=\phi(u), \ \ \alpha a+\beta=b
\end{equation*}

It is well known that if a function with an algebraic addition theorem satisfies an equation $\phi(\alpha u+\beta)=\phi(u)$, then $\alpha$ is a root of unity of arbitrarily high exponent for an algebraic function, can only have the values $\pm 1$ for a simply periodic function, and for a doubly periodic function can only have one of the values $\pm 1$, $\pm i$, $\pm\frac{1}{2}\pm\frac{1}{2}\sqrt{3}i$.\footnote{For this see my paper ``Ueber die Uniformisierung der algebraischen Kurven 1".  Math. Annalen, Vol 67 (1909), in particular p. 164 ff.}

%%%%%%%%%%%%%%%%%%%%%%%%%%%%%%%%%%%%%%%%%%%%%%%%

\newpage

\section{The degree of the irreducible equation \newline $G[\phi(u),  \phi(v),\phi(u+v)]=0$ with respect to $\phi(u+v)$, $\phi(u)$, $\phi(v)$}

Again, we let $\phi(u)$ be a function with an algebraic addition theorem.  A simple formula can be found for the degree of the irreducible algebraic equation which expresses the addition theorem with respect to $\phi(u+v)$, $\phi(u)$, $\phi(v)$.

Let $\nu$ be the order of the function  $\phi(u)$.  Then to an arbitrarily chosen pair of values  $\phi(u)$,  $\phi(v)$, there belong $\nu$ incongruent values of $u$ and $\nu$ incongruent values of $v$.  Accordingly we obtain $\nu^2$ incongruent values of $u+v$ and, if $m$ is the number of branches of the function  $\phi(u)$, a total of $m\nu^2$ values of  $\phi(u+v)$.  By an application of a method given in \S 3 it can be shown that the known $m\nu^2$ values  $\phi(u+v)$ are the values of just as many branches of one and the same algebraic function of the independent variables  $\phi(u)$ and  $\phi(v)$.  Accordingly, this algebraic function possesses c different branches if no two of the known $m\nu^2$ values of $\phi(u+v)$ are equal to each other for all values of $\phi(u)$ and $\phi(v)$.  The algebraic function under consideration cannot have more than $m\nu^2$  values.  Otherwise it can happen in particular cases that two of the $m\nu^2$ values $\phi(u+v)$ ar equal to each other for all values of $\phi(u)$ and $\phi(v)$.  This can happen if and only if it is possible to determine four variables $u_1$, $v_1$, $u_2$, $v_2$, of which only the first two are independent, such that the following equations hold\begin{equation*}
\label{ }
\phi(u_1)=\phi(u_2), \ \ \phi(v_1)=\phi(v_2), \ \ \phi(u_1+v_1)=\phi(u_2+v_2)
\end{equation*}without it being true that simultaneously $u_1$ is congruent to $u_2$ and $v_1$ is congruent to $v_2$.  

If we apply the results in \S 4 for the equations $$\phi(u)=\psi(u_1), \ \ \phi(v)=\psi(v_1), \ \ \phi(u+v)=\psi(u_1+v_1)$$ to the previous equations above , it must be true that \begin{equation*}
\label{ }
u_1\ \  \text{congruent}\ \ \alpha u_2 \ \ \text{and} \ \ v_1\ \  \text{congruent}\ \  \alpha v_2
\end{equation*}where $\alpha$ is a constant, i.e., if two of the $m\nu^2$ values of $\phi(u+v)$ are equal for all values of $\phi(u)$ and $\phi(v)$ the function has the property of remaining unchanged after multiplication of its argument by a constant different from \emph{one}.

For an arbitrary function $\phi(u)$ with an algebraic addition theorem we denote the number of different constant factors $\alpha$, (including $\alpha=1$) for which $\phi(\alpha u)=\phi(u)$ by $\lambda_0$.  According to the remarks made at the end of \S 4, the number $\lambda_0$ can have any positive integral value for an algebraic function, only one of the values $1,2$ for a simply periodic function, and only one of the values $1,2,3,4,6$ for a doubly periodic function.  In every case, $\dfrac{\nu}{\lambda_0}$ is an integer.  It follows that the set of  $m\nu^2$ values of $\phi(u+v)$ are broken into $m \nu\cdot \dfrac{\nu}{\lambda_0} $ classes with $\lambda_0$ members falling into each class, and to be sure in such a way that for any two values $\phi(u+v)$ which belong to different classes do not have all of the values of $\phi(u)$ and $\phi(v)$ equal to each other.  Therefore, it has been shown:

\emph{The degree of the irreducible equation $G[\phi(u),  \phi(v),\phi(u+v)]=0$ with respect to $\phi(u+v)$ is equal to $m \nu\cdot \dfrac{\nu}{\lambda_0} $.}

The degree of the equation $G[\phi(u),  \phi(v),\phi(u+v)]=0$ with respect to $\phi(u)$ or to $\phi(u)$ is also equal to $m \nu\cdot \dfrac{\nu}{\lambda_0} $, as can be proven by a process analogous to the preceding.

The theorem just proven leads to the solution of the following problem:  To determine all analytic functions of the argument $u$ for which $\phi(u+v)$ is rationally expressible through $\phi(u)$ and $\phi(v)$ alone.\footnote{The same problem, which I will subsequently solve, was solved by \textsc{Weierstrass} in another way in his lectures on the theory of elliptic functions held in the winter semester 1862/63.}  For these functions must $m \nu\cdot \dfrac{\nu}{\lambda_0} =1$, and consequently $m=\nu=1$, whence:

\emph{If an analytic function  $\phi(u)$ has the property that for it $\phi(u+v)$ is rationally expressible by $\phi(u)$ and $\phi(v)$ alone, then $\phi(u)$ is either a linear or linear fractional function of $u$ or a linear or linear fractional function of $e^{\mu u}$}.
\\
\\
\hrule
\bigskip
\bigskip

While the irreducible algebraic equation $(G)$ relating the function values $\phi(u)$, $\phi(v)$, $\phi(u+v)$ will generally change form if we move the origin of the complex plane, the irreducible algebraic equation \begin{equation*}
\label{ }
K[\phi(u),\phi(v);\phi(w),\phi(t)]=0
\end{equation*}where the relation $$u+v=w+t$$holds, obtained from the corresponding equation from $(G)$ by elimination, does not change its form if we move the origin.  By means of reasoning analogous to that which we used for the equation $G[\phi(u),  \phi(v),\phi(u+v)]=0$ we obtain:

\emph{The degree of the irreducible equation\begin{equation*}
\label{ }
K[\phi(u),\phi(v);\phi(w),\phi(t)]=0
\end{equation*}with respect to each of the quantities  $\phi(u)$, $\phi(v)$, $\phi(w)$, $\phi(t)$ is the same, namely $m\dfrac{\nu^2}{\lambda}$.}

The number $\lambda$, which is a divisor of $\nu$, has the following meaning.  To the function $\phi(u)$ belongs a definite \emph{group} of linear substitutions $u'=\alpha u+\beta$ of the argument $u$, consisting of all those substitutions of the given form for which $\phi(u')$=$\phi(u)$.

The different coefficients $\alpha$ occuring in this group of substitutions are finite in number and are all roots of unity.  If $\epsilon$ is one of these roots of unity and, to be sure, one of the highest possible order, this order is equal to $\lambda$ and the set of coefficients $\alpha$ coincides with the set of $\lambda^{\text{th}}$ roots of unity.\footnote{see the citation on page 12.}
\\
\hrule\hrule

\end{document}